\magnification=1100
\font\bigbf=cmbx10 scaled \magstep2
\font\medbf=cmbx10 scaled \magstep1 
\hfuzz=10 pt
\input psfig.sty

\def\i{\item}
\def\n{\noindent}
\def\L{{\bf L}}

\def\n{\noindent}

\def\Ll{{\bf L}^1_{\rm loc}}

\def\be{{\bf e}}

\def\c{\centerline}

\def\R{I\!\!R}

\def\Tilde{\widetilde}

\def\vs{\vskip 2em}
\def\vsk{\vskip 4em}
\def\v{\vskip 1em}
\null
\c{\bigbf An Ill Posed  Cauchy Problem}
\v
\c{\bigbf for a  Hyperbolic System in Two Space Dimensions}
\vs
\c{\it Alberto Bressan}
\v
\centerline{S.I.S.S.A. - Via Beirut 4, Trieste 34014, Italy}
\centerline{e-mail: bressan@sissa.it}
\vs
\n{\medbf 1 - Introduction}
\v
The theory of weak solutions for nonlinear conservation laws is 
now well developed in the case of scalar equations [3] and for
one-dimensional hyperbolic systems [1, 2].   For systems in several space
dimensions, however, even the global existence of solutions to the 
Cauchy problem remains a challenging open question.
In this note we construct a conterexample showing that, even for
a simple class of hyperbolic systems, in two space dimensions
the Cauchy problem can be ill posed.
\v
The systems on $\R^m$ that we consider take
the special form
$${\partial\over\partial t} u_i+\sum_{\alpha=1}^m{\partial\over
\partial x_\alpha}
\Big(f_\alpha\big(|u|\big)\,
u_i\Big)=0\qquad\qquad  i=1,\ldots,n\,.\eqno(1.1)$$
For physical motivations, see [2] or [5].
Given a sufficiently regular initial condition
$$u(0,x)=\bar u(x)\,,\eqno(1.2)$$
the solution of the Cauchy problem can be constructed as follows.
\v
\n{\bf 1.} Set $\rho\doteq |u|$ and solve the Cauchy problem
for a scalar conservation law
$$\rho_t+\sum_{\alpha=1}^m\big(f_\alpha(\rho)\,\rho\big)_{x_\alpha}=0\,,
\qquad\qquad \rho(0,x)=\big|\bar u(x)\big|\,.\eqno(1.3)$$
This will provide the absolute value of the solution.
\v
\n{\bf 2.} To find the angular component $\theta\doteq u/|u|$,
consider the O.D.E.
$$\dot x=f\big(\rho(t,x)\big)\eqno(1.4)$$
with $f=(f_1,\ldots,f_m)$.
Assuming that (1.4) has a unique solution for almost every initial 
data $x(0)=y$, 
denote by $t\mapsto x(t)\doteq \Phi_t y$ the corresponding trajectory.
\v
\n{\bf 3.} Given $t>0$ and $x\in\R^m$, call 
$y\doteq \Phi^{-t}x$ be the (unique) point such that
$\Phi_t y=x$.  Assuming that the measurable maps $\Phi^{-t}$
are well defined (up to sets of zero measure),
the angular component of the solution of (1.1)-(1.2) is then obtained as
$$\theta(t,x)=\bar \theta\big(\Phi^{-t}x\big)\,,\eqno(1.5)$$
where $\bar\theta\doteq \bar u/|\bar u|$.
\vs
If the solution $\rho=\rho(t,x)$ of 
the scalar conservation law
(1.3) remains  smooth (or piecewise smooth), 
it is easy to check that the above 
procedure actually yields a weak solution $u$ to the Cauchy problem 
(1.1)-(1.2).  
In general however, even if the initial data $\rho(0,\cdot)$
is smooth, regularity may be lost in finite time.
It is thus interesting to understand what kind of
assumptions can provide the existence of entropy weak solutions.
A natural set of conditions is the following:
\v
\i{(A1)} The flux function $F(\rho)\doteq f(\rho)\,\rho$
is Lipschitz continuous.
\v
\i{(A2)} The initial data is bounded, measurable, and bounded 
away from zero, i.e.
$$\bar u\in\L^\infty\,,\qquad\qquad 
0<a\leq\big|\bar u(x)\big|\leq b\qquad\hbox{for a.e.~}x\in\R^m\,.$$   
\v
In the case of one space dimension,
the above assumptions yield the global existence of
a unique entropy weak solution to the Cauchy problem.
Indeed, by the fundamental theorem of Kruzhkov [3],
the equation (1.3) admits a unique entropy weak solution
$\rho=\rho(t,x)$, satisfying
$$a~\leq~\rho(t,x)~\leq ~b\,.\eqno(1.6)$$
Let $x_1(t)< x_2(t)$
be any two solutions of (1.4).  For every $t>0$
the conservation equation implies
$$\int_{x_1(t)}^{x_2(t)}\rho(t,x)\,dx~=~
\int_{x_1(0)}^{x_2(0)}\rho(0,x)\,dx\,.$$
By (1.6) this yields
$${a\over b}\big(x_2(0)-x_1(0)\big)~\leq~ x_2(t)-x_1(t)~\leq~
{b\over a}\big(x_2(0)-x_1(0)\big).$$
Hence the maps 
$$(t,y)\mapsto \Phi_t y\,,\qquad\qquad (t,x)\mapsto
\Phi^{-t}x$$
are both well defined and Lipschitz continuous.
It is now clear that (1.5) yields the desired solution.
\v
The example constructed below will show that, under exactly the same 
assumptions (A1)-(A2) which provide the  existence and uniqueness
of solutions in the one-dimensional case, 
in two space dimensions the Cauchy problem (1.1)-(1.2)
is not well posed. 
The reason is that, although the scalar conservation law (1.3)
always admits a unique entropy solution,
the corresponding O.D.E.~in (1.4) has discontinuous right hand side.
For a given time $t>0$, the flux function $y\mapsto\Phi_t y$
may not admit a well defined inverse,
because of oscillation phenomena.
\vsk
\n{\medbf 2 - A counterexample}
\v
Call $\be_1 =(1,0)$, $\be_2=(0,1)$
the canonical basis on $\R^2$.
Consider the Lipschitz (piecewise affine) flux function
$$F(\rho)\doteq \cases{
0\qquad &if\quad $\rho\leq 1\,$,\cr
(1-\rho)\be_1\qquad &if\quad $1\leq\rho\leq 2\,$,\cr
(\rho-3)\be_1\qquad &if\quad $2\leq \rho\leq 3\,$,\cr
(\rho-3)\be_2\qquad &if\quad $3\leq \rho\,$.\cr}\eqno(2.1)$$
Observe that
$$F'(\rho)= \cases{
0\qquad &if\quad $\rho\leq 1\,$,\cr
-\be_1\qquad &if\quad $1\leq\rho\leq 2\,$,\cr
\be_1\qquad &if\quad $2\leq \rho\leq 3\,$,\cr
\be_2\qquad &if\quad $3\leq \rho\,$.\cr}\eqno(2.2)$$
We can write $F(\rho)=f(\rho)\,\rho$, with
$$f(\rho)\doteq 
\cases{0\qquad &if\quad $\rho\leq 1\,$,\cr
\rho^{-1}(1-\rho)\be_1\qquad &if\quad $1\leq\rho\leq 2\,$,\cr
\rho^{-1}(\rho-3)\be_1\qquad &if\quad $2\leq \rho\leq 3\,$,\cr
\rho^{-1}(\rho-3)\be_2\qquad &if\quad $3\leq \rho\,$.\cr}\eqno(2.3)$$
In particular,
$$f(2)=-{1\over 2}\be_1\,,\qquad\qquad f(3)=0\,,\qquad\qquad f(4)=
{1\over 4}\be_2\,.\eqno(2.4)$$
We now look at the O.D.E.
$$\dot x=f\big(\rho(t,x)\big)\eqno(2.5)$$
in connection with
special types of solutions $\rho=\rho(t,x)$ of the conservation law
$$\rho_t+\nabla_x\cdot F(\rho)=0\,.\eqno(2.6)$$
For a given set $Q\subset \R^2$, 
we consider three entropy weak solutions of (2.6), namely 
$$\rho^\natural(t,x)\doteq 3\,,$$

$$\rho^\sharp(t,x)\doteq \cases{4\qquad &if\quad $x-t\be_2\in Q\,$,\cr
3\qquad &otherwise,\cr}\eqno(2.7)$$

$$\rho^\flat(t,x)\doteq \cases{2\qquad &if\quad $x-t\be_1\in Q\,$,\cr
3\qquad &otherwise.\cr}\eqno(2.8)$$

\midinsert
\vskip 10pt
\centerline{\hbox{\psfig{figure=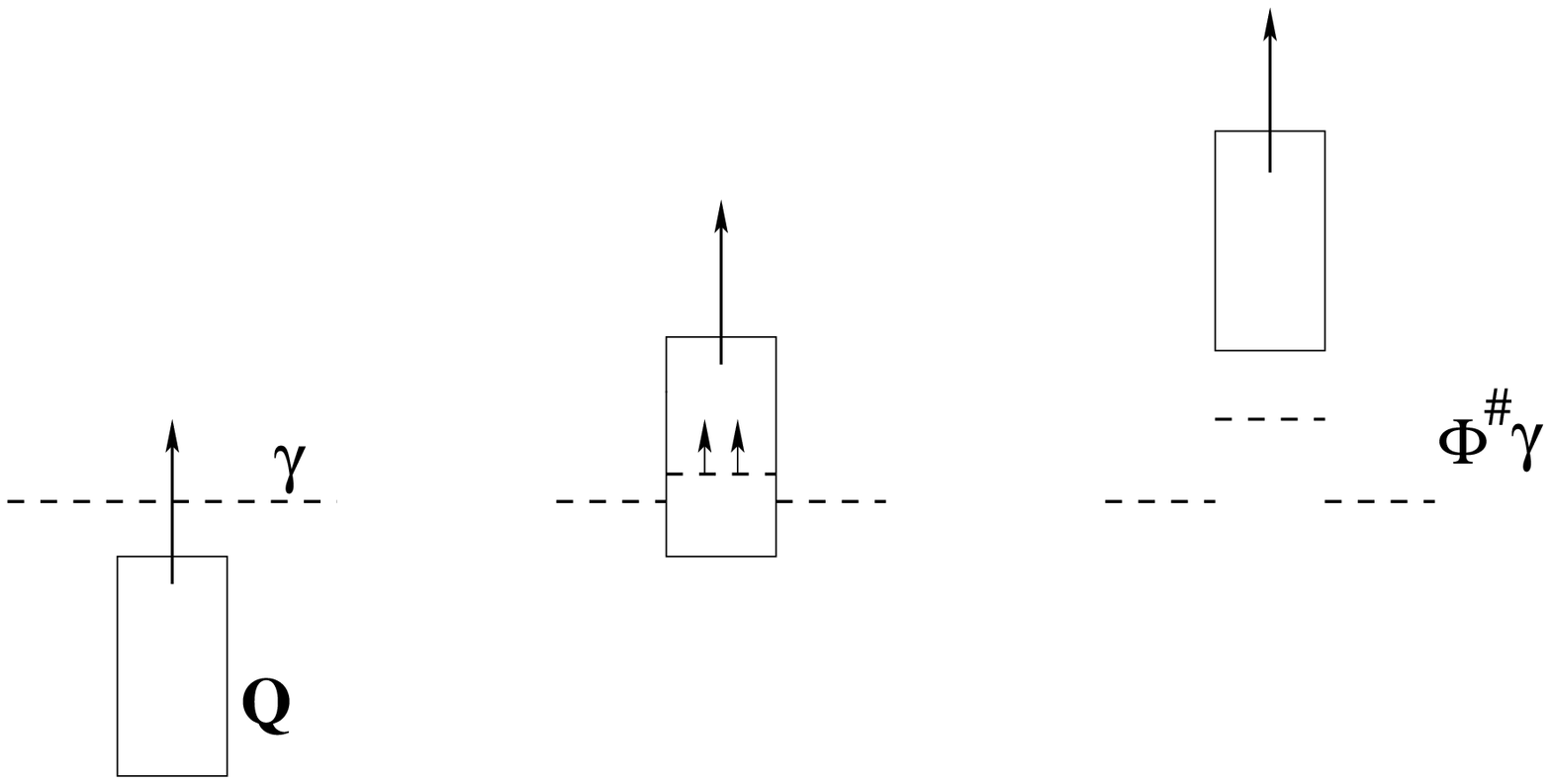,width=8cm}}}
\centerline{\hbox{figure 1}}
\vskip 10pt
\endinsert

We begin by studying what happens to the corresponding
trajectories of (2.5)
in the special case where $Q$ is a rectangle:
$$Q\doteq [0,a]\times [0,b]\subset \R^2\,.$$

\n{\bf 1.} In the case $\rho=\rho^\natural$, we trivially have
$$\dot x=0\qquad\qquad\Phi^\natural_t y=y\,.$$
\v
\n{\bf 2.} When
$\rho=\rho^\sharp$, we have (fig.~1)
$$\dot x=\cases{\be_2/4\qquad &if\quad $x_1\in [0,a]$, 
~~$x_2-t\in [0,b]$,\cr
0\qquad &otherwise.\cr}$$
Points on the vertical strip 
$$\Gamma_a\doteq\big\{ (x_1,x_2)\,;~~x_1\in [0,a]\big\}
$$
are eventually
shifted to the right 
by an amount $b\be_1/3$, while all other points do not move.
More precisely:
$$\Phi^\sharp_t y-y =\cases{b\be_1/3\qquad &if\qquad $y_1\in [0,a],
~~y_2\geq b,~~t\geq (y_2-b)+4b/3\,,$\cr
0\qquad &otherwise.\cr}$$
This is illustrated in fig.~1.  
One has $\rho=4$ on $Q$ and $\rho=3$ outside.
The rectangle $Q$ moves
upward with unit speed.  At any given time $t$, points which lie inside the
rectangle move with speed $\dot x=\be_2/4$. 
The points which initially lie on the segment $\gamma$ 
are eventually displaced onto the set $\Phi^\sharp \gamma$.
\v
\n{\bf 3.}
When $\rho=\rho^\flat$, we have
$$\dot x=\cases{-\be_1/2\qquad &if\quad $x_1-t\in [0,a]$,
~~$x_2\in [0,b]$,\cr
0\qquad &otherwise.\cr}$$
Points on the horizontal strip 
$$\Gamma_b\doteq\big\{ (x_1,x_2)\,;~~x_2\in [0,b]\big\}
$$
are eventually
shifted to the left 
by an amount $-a\be_1/3$, while all other points do not move.
More precisely:
$$\Phi^\flat_t y-y =\cases{-a\be_1/3\qquad &if\qquad $y_2\in [0,b],
~~y_1\geq a,~~t\geq (y_1-a)+2a/3\,,$\cr
0\qquad &otherwise.\cr}$$
\v
Starting with these elementary solutions, we can construct more 
complex ones.
For example, in (2.7)  can take
$$Q=Q^\sharp\doteq \bigcup_{k~\rm{even}} [k,k+1]\times [k,\,k+3]$$
In the corresponding equation (2.5), this generates the shift
$$\Phi^\sharp y-y=\cases{ \be_2\qquad &if\quad $[[y_1]]$ is even,\cr
0\qquad &if\quad $[[y_1]]$ is odd.\cr}$$

On the other hand, in (2.8) we can take
$$Q=Q^\flat\doteq \bigcup_{k~\rm{even}} [k,k+3]\times [k,\,k+1]$$
In the corresponding equation (2.5), this generates the shift
$$\Phi^\flat y-y=\cases{ -\be_1\qquad &if\quad $[[y_2]]$ is even,\cr
0\qquad &if\quad $[[y_2]]$ is odd,\cr}$$
where $[[s]]$ denotes the largest integer $\leq s$.

\midinsert
\vskip 10pt
\centerline{\hbox{\psfig{figure=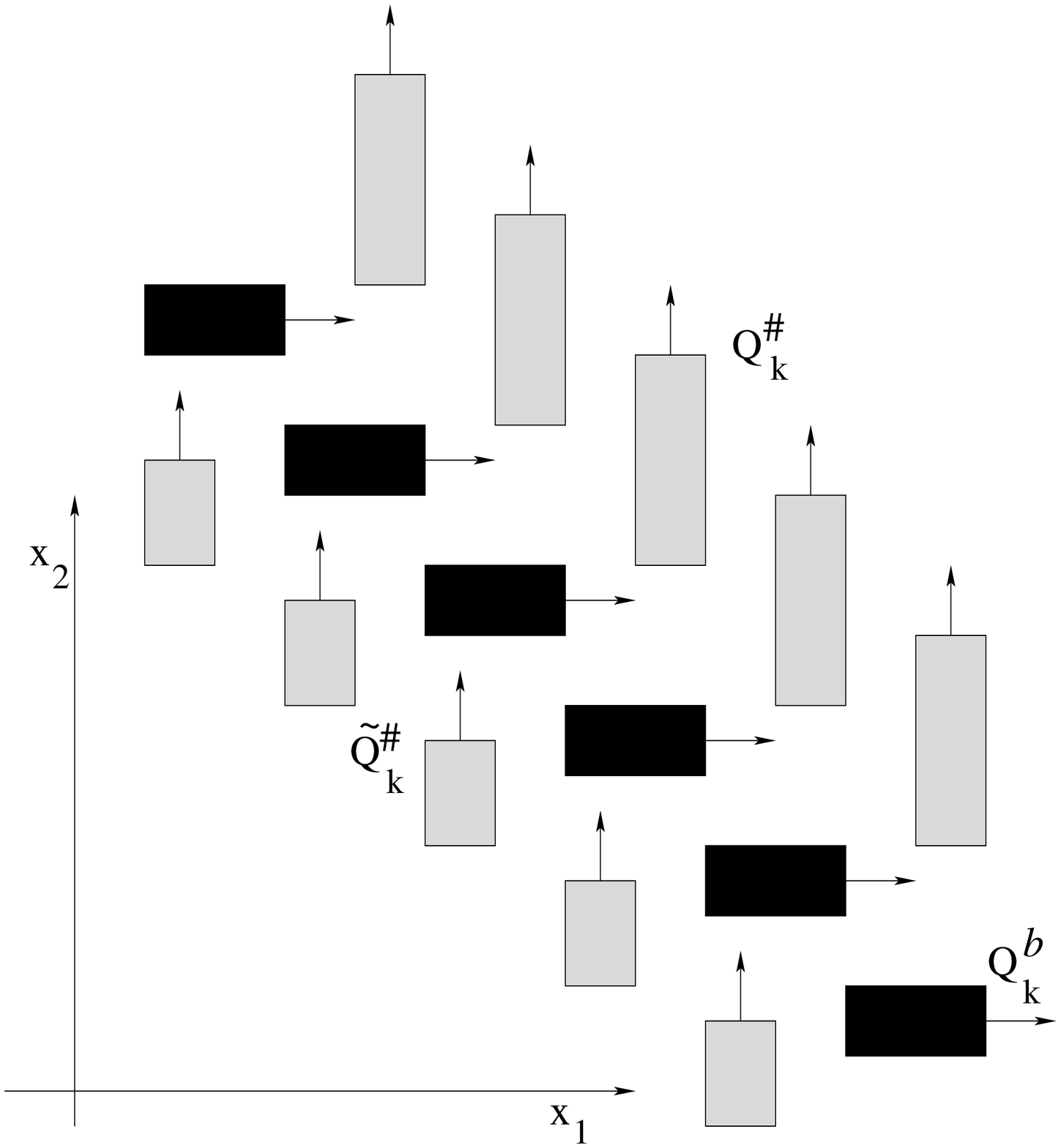,width=10cm}}}
\centerline{\hbox{figure 2}}
\vskip 10pt
\endinsert

In the following, we write the coordinates of a point $x=(x_1,x_2)\in\R^2$
in binary digits, say
$$x_1=\sum_k \alpha_k 2^{-k}\,,\qquad\qquad 
x_2=\sum_k \beta_k 2^{-k}\,.\eqno(2.9)$$
with $\alpha_k,\beta_k\in\{0,1\}$.
For every $k$, the previous analysis shows that one can construct
sets $Q^\sharp_k$, $Q^\flat_k$, $\Tilde Q^\sharp_k$ 
contained inside disjoint strips
$$\eqalign{
Q^\sharp_k&\subseteq \big\{ (x_1,x_2)\,;~~x_2+x_1\in 2^{-k}\,[28,\,32]\big\}\,,
\cr
&\cr
Q^\flat_k&\subseteq \big\{ (x_1,x_2)\,;~~x_2+x_1\in 2^{-k}\,[22,\,27]\big\}\,,
\cr
&\cr
\Tilde 
Q^\sharp_k&\subseteq \big\{ (x_1,x_2)\,;~~x_2+x_1\in 2^{-k}\,[16,\,21]
\big\}\,,
\cr}\eqno(2.10)$$
and such that the following holds.
\v
\n{\bf (i)}~~Given the solution
$$\rho_k^\sharp(t,x)\doteq \cases{4\qquad &if\quad $x-t\be_2\in Q^\sharp_k
\,$,\cr
3\qquad &otherwise,\cr}$$
the corresponding equation (2.5) generates the shift 
$$\Phi^\sharp y-y=\cases{ 2^{-k}\be_2\qquad &if\qquad $\alpha_k=0$,\cr
0\qquad &if\qquad $\alpha_k=1$.\cr}$$
\v
\n{\bf (ii)}~~Given the solution
$$\rho_k^\flat(t,x)\doteq 
\cases{2\qquad &if\quad $x+t\be_1\in Q^\flat_k\,$,\cr
3\qquad &otherwise,\cr}$$
the corresponding equation (2.5) generates the shift 
$$\Phi^\flat y-y=\cases{ -2^{-k}\be_1\qquad &if\qquad $\beta_k\not=
\beta_{k+1}$,\cr
0\qquad &if\qquad $\beta_k=\beta_{k+1}$.\cr}$$
\v
\n{\bf (iii)}~~Given the solution
$$\tilde
\rho_k^\sharp(t,x)\doteq \cases{4\qquad &if\quad $x-t\be_2\in \Tilde
Q^\sharp_k\,$,\cr
3\qquad &otherwise,\cr}$$
the corresponding equation (2.5) generates the shift 
$$\Tilde\Phi^\sharp y-y=\cases{ 2^{-k-1}
\be_2\qquad &if\qquad $\alpha_k=0$,\cr
0\qquad &if\qquad $\alpha_k=1$.\cr}$$
\v
As shown in fig.~2, the set $Q_k^\sharp$
consists of an array of (gray) 
rectangles of size $2^{-k}\times 3\cdot 2^{-k}$,
the set $Q_k^\flat$ consists of an array of (black) 
rectangles of size $2^{1-k}\times 2^{-k}$, while
the set
$\Tilde Q_k^\sharp$ consists of an array of (gray) 
rectangles of size $2^{-k}\times 3\cdot 2^{1-k}$.
Moving the gray rectangles upward with unit speed and moving
the black ones
to the right with unit speed, the three arrays will never overlap.
We can thus consider the composite solution
$$\rho_k(t,x)\doteq \cases{4\qquad &if\quad $x-t\be_2\in Q^\sharp_k\cup
\Tilde  Q^\sharp_k\,$,\cr
2\qquad &if\quad $x-t\be_1\in Q^\flat_k\,$,\cr
3\qquad &otherwise.\cr}\eqno(2.11)$$
The differential equation
$$\dot x=f\big(\rho_k(t,x)\big)$$
will eventually shift points on the plane, according to
the composition
$$\Psi_k\doteq\Tilde\Phi^\sharp_k\circ\Phi^\flat_k\circ\Phi^\sharp_k\,.
\eqno(2.12)$$

\midinsert
\vskip 10pt
\centerline{\hbox{\psfig{figure=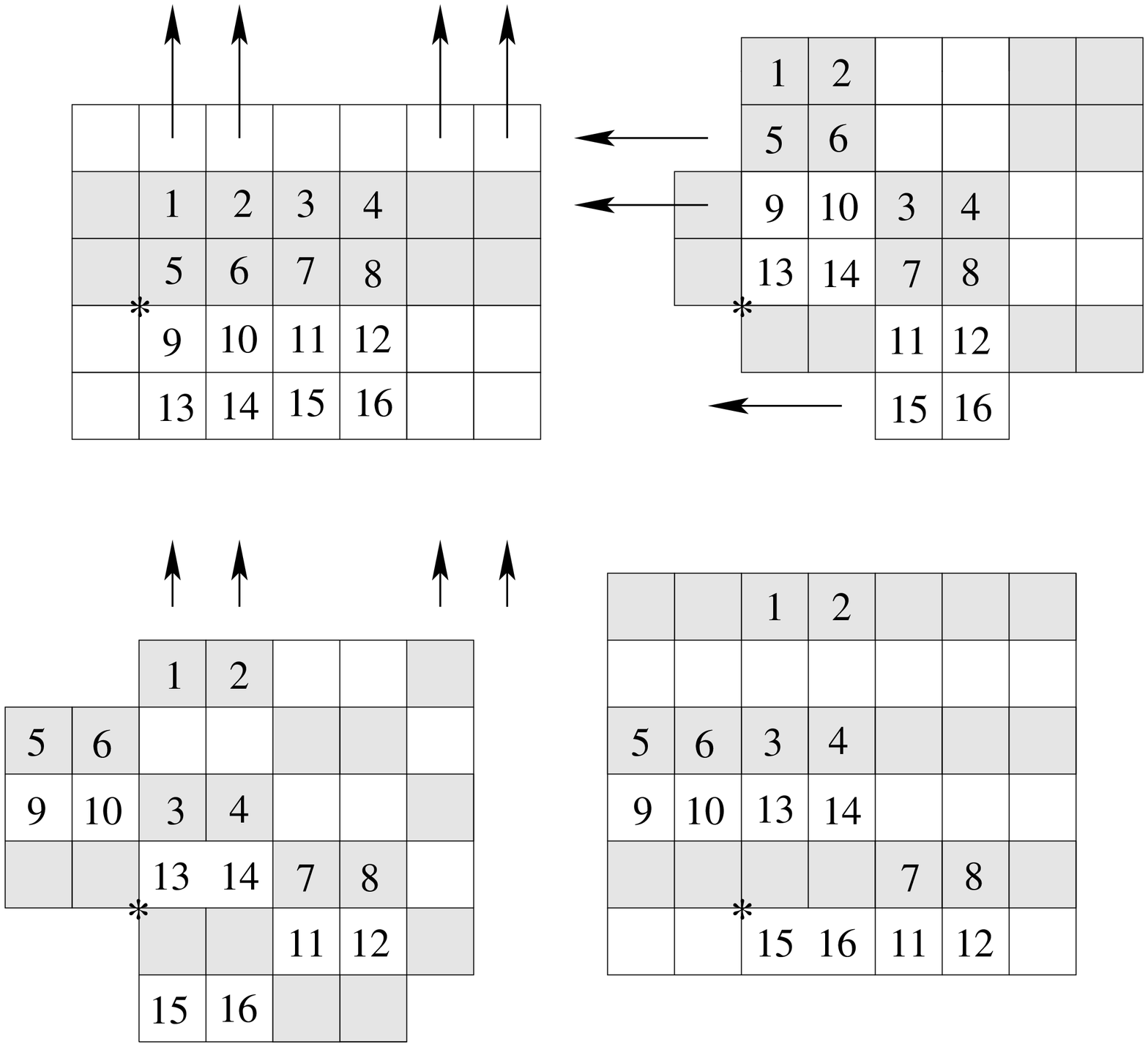,width=10cm}}}
\centerline{\hbox{figure 3}}
\vskip 10pt
\endinsert

This is illustrated in fig.~3. Each square here has sides of length
$2^{-k-1}$. In all four pictures, the asterisque marks the samepoint 
in the plane.
Performing first the upward shift $\Phi^\sharp_k\,$, then the leftward
shift $\Phi_k^\flat$ and finally the upward shift 
$\Tilde\Phi^\sharp_k\,$,
one obtains a measure-preserving transformation $\Psi_k$
whose main property is the following.
In connection with the dyadic decomposition (2.9), define
$$X_j\doteq \big\{ (x_1,x_2)\,;~~\beta_j=0\big\},\qquad\qquad 
X'_j\doteq \big\{ (x_1,x_2)\,;~~\beta_j=1\big\}.$$
Then, up to sets of measure zero, one has
$$\Psi_k(X_k)=X_{k+1}\,,\qquad\qquad \Psi_k(X'_k)=X'_{k+1}\,.\eqno(2.13)$$
Intuitively, one can think at the moving sets $Q^\sharp_k$, $Q^\flat_k$,
$\Tilde Q^\sharp_k$ as a comb, moving in the north-east direction.
The effect of its passage is to displace points on the set
$X_k$ onto the set $X_{k+1}$.
Putting an arbitrarily large number of these combs one next to the other
(fig.~4), 
we now construct a sequence of solutions $u_n$
where the initial data converge strongly in $\Ll$ but at a later time
$\tau$ the corresponding solutions have no strong limit.
Hence the limiting Cauchy problem
cannot be well posed.
\v
Define a solution of (1.3) by setting
$$\rho_n(t,x)\doteq \cases{4\qquad &if\quad $x-t\be_2\in \bigcup_{0\leq k\leq
n} 
\big(Q^\sharp_k\cup
\Tilde  Q^\sharp_k\big)\,$,\cr
&\cr
2\qquad &if\quad $x-t\be_1\in  \bigcup_{0\leq k\leq
n} Q^\flat_k\,$,\cr
&\cr
3\qquad &otherwise.\cr}\eqno(2.14)$$
Consider the Cauchy problem (1.1)-(1.2) with
$f$ given at (2.3).   We choose an initial value $\bar u_n$ of the form
$$\bar u_n(x)=\rho_n(0,x)\cdot
\big(\cos \bar\theta(x)\,,~\sin\bar\theta(x)\big)\eqno(2.15)$$
with 
$$\bar\theta(x)\doteq\cases{\beta\qquad &if\quad $[[x_2]]$ is even,\cr
-\beta\qquad &if\quad $[[x_2]]$ is odd,\cr}$$
for a fixed angle $\beta\in\,]0,\pi/2[\,$.
We claim that the sequence of weak solutions $u_n$ has no limit in 
$\Ll$ as $n\to\infty$. 
Indeed, for $32<x_1+x_2< t$
one has
$$u_n(t,x)=3\big(\cos\theta_n(t,x)\,,~\sin\theta_n(t,x)\big),$$
with $x=(x_1,x_2)$ and
$$\theta_n(t,x)=\cases{\beta\qquad &if\quad $[[2^n x_2]]$ is even,\cr
-\beta\qquad &if\quad $[[2^n x_2]]$ is odd.\cr}$$
As a result, when $n\to\infty$, the 
$u_n$ do not converge in $\Ll$ and the
Cauchy problem is not well posed.
\v

\midinsert
\vskip 10pt
\centerline{\hbox{\psfig{figure=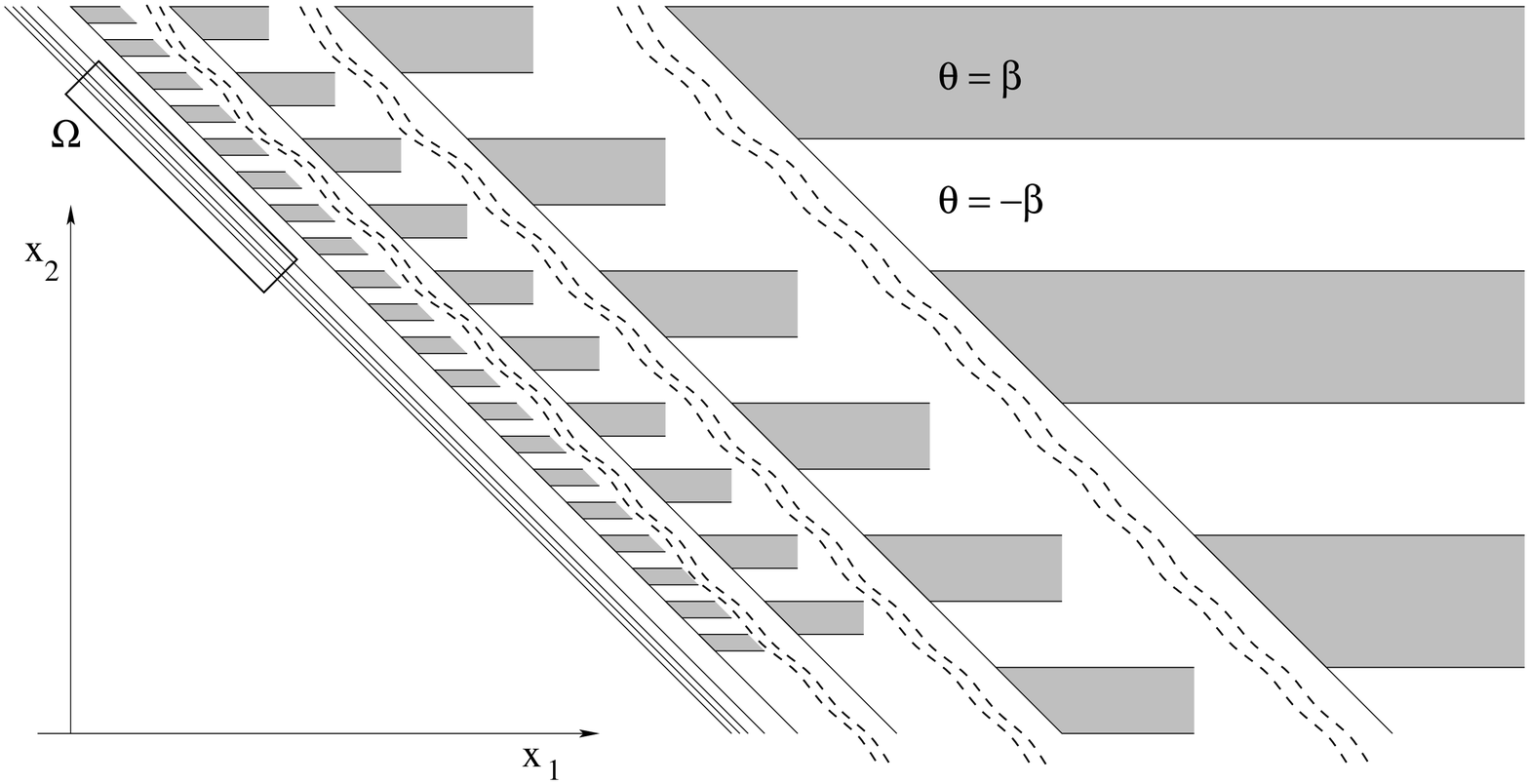,width=12cm}}}
\centerline{\hbox{figure 4}}
\vskip 10pt
\endinsert

On the other hand,
still in the range where $32<x_1+x_2< t$, one can compute 
the weak limit of the sequence $\{u_n\}$, given by
$$u(t,x)=(3\cos\beta\,,~0).$$
To check whether this weak limit provides a distributions solution
to the system (1.1), (2.3), it suffices to examine what happens
in a neighborhood of the line $x_1+x_2=t$.
Computing the time derivative
$${d\over dt}\int_{\Omega(t)} u(t,x)\,dx$$
for a set of the form
$$\Omega(t)\doteq\big\{ (x_1,x_2)\,,~~~x_2\in [a,b]\,,~~x_1\in [t-x_2-2^{-k}\,,
t-x_2+2^{-k}]\big\},$$
one checks that the conservation equations are satisfied
provided that
$F(3\cos\beta)=F(3)=0$. In other words:
\v
\n - If $\cos \beta\leq 1/3$, then the weak limit $u$ provides a 
weak solution to the system of conservation laws.   
This yields a curious example where entropy is dissipated by 
a linearly degenerate field.
\v
\n - If $\cos \beta> 1/3$, then the weak limit $u$ is not
weak solution.   In this case, it is not even clear whether any
solution exists at all.
\vsk
\n{\medbf 3 - Concluding remarks}
\v
\n{\bf 1.} All systems of the form (1.1) are hyperbolic,
but in a certain way pathological.
In one space dimension, the system
$$u_t+\Big(f\big(|u|\big)u\Big)_x=0\eqno(3.1)$$
takes the quasilinear form
$$u_t+A(u)u_x=0,\qquad\qquad A(u)=f\big(|u|\big)\, I+f'\big(|u|\big)\,
u\otimes u\,.$$
At a given point $u\not= 0$, the Jacobian matrix $A(u)$
has the 1-dimensional eigenspace
$E\doteq \R u$
spanned by the vector $u$, corresponding to the eigenvalue 
$\lambda=f\big(|u|\big)+f'\big(|u|\big)\,|u|$.
Moreover, one finds a second eigenvalue
$\lambda^*\doteq f\big(|u|\big)$
of multiplicity $n-1$, corresponding to
the orthogonal eigenspace
$E^\perp$ (fig.~5).
\midinsert
\vskip 10pt
\centerline{\hbox{\psfig{figure=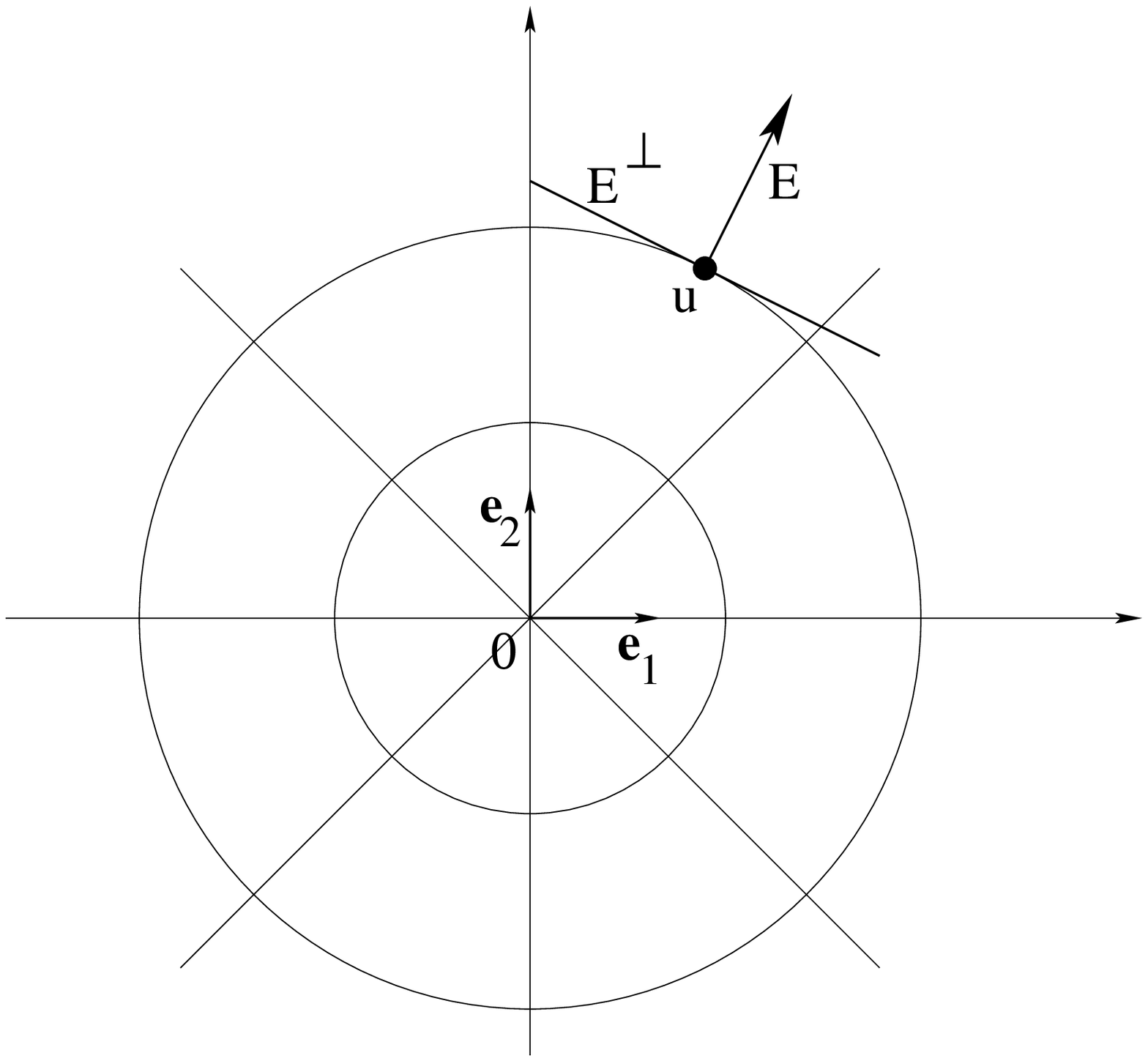,width=8cm}}}
\centerline{\hbox{figure 5}}
\vskip 10pt
\endinsert
Even in one space dimension
the Riemann problem can have multiple solutions, 
all obtained as vanishing viscosity limits.
For example, consider
$$u_t+\big(|u^2|u\big)_x=0\qquad\qquad u\in\R^2,\eqno(3.2)$$
$$u(0,x)=\cases{\be_1\qquad &if\quad $x<0$,\cr
-\be_1\qquad &if\quad $x>0$.\cr}$$
An entropic solution is provided by
$$u(t,x)=\rho(t,x)\be_1\,,$$
where $\rho$ is the entropy weak solution of the corresponding scalar problem
$$\rho_t+(\rho^3)_x=0\qquad\qquad 
\rho(0,x)=\cases{1\qquad &if\quad $x<0$,\cr
-1\qquad &if\quad $x>0$.\cr}$$
This is clearly a limit of vanishing viscosity approximations [3].

In addition, there is the second admissible solution
$$u(t,x)=\cases{\be_1\qquad &if\quad $x<t$,\cr
-\be_1\qquad &if\quad $x>t$.\cr}\eqno(3.3)$$
Indeed, consider any smooth function
$\theta:\R\mapsto [0,\pi]$ such that 
$$\theta(s)=\cases{0\quad &if\quad $x\leq 0$,\cr
\pi\qquad &if\quad $x\geq 1$.\cr}$$
Then the Cauchy problem for (3.2) with initial data
$$u_n(0,x)=\big(\cos \theta(nx)\,,~\sin\theta(nx)\big)$$
has a unique smooth solution, namely
$u_n(t,x)~=~u(0, x-t)$.  This is also a limit of vanishing viscosity
approximations.  Letting $n\to\infty$, these solutions 
converge to (3.3).  

We remark, however, that our counterexample 
is not related to the presence of a singularity at the origin.
Indeed, the initial data at (2.14)-(2.15) 
are contained in a bounded domain whose convex
closure does not contain the origin.
\vs
\n{\bf 2.} The flux function is not smooth but only Lipschitz continuous.
The oscillations in the solution of the scalar conservation law
(1.3) are not damped, because of the linear degeneracy of the
flux.  This indeed leaves open the possibility that, with a smooth flux,
such a pathological behavior will not be observed.

We recall that, in the one dimensional case, it is common to
replace a 
general flux $F$ by a polygonal approximation $F_\nu$. 
This is indeed the basis for the method 
of front-tracking approximations.  
In one space dimension, the Cauchy problems related to 
all polygonal fluxes are well posed [1].   Taking a limit,
one thus obtains 
a proof of well posedness for the original Cauchy problem.
The multidimensional case is thus completely different.
\vs
\n{\bf 3.} What our example points out is that, to get well posedness,
the initial data should be chosen in a space of functions
smaller than $\L^\infty$.
In view of the one-dimensional theory, it is natural
to conjecture that a global weak solution
should exist for initial data $\bar u$ with bounded variation.
A proof of this conjecture (far from easy!) could rely
on a compactness property for fluxes generated by
O.D.E's with right hand side in BV.  More precisely,
consider a sequence of smooth maps $f_\nu:[0,T]\times\R^m\mapsto\R^m$
such that
$$\big|f_\nu(t,x)\big|\leq C_1\,,\eqno(3.4)$$
$$\|f_\nu\|_{\strut BV}~\doteq~\int_0^T\int_{\R^m} 
\left|{\partial\over\partial t}f_\nu\right|+\sum_{i=1}^m
\left|{\partial\over\partial x_i}f_\nu\right|~~dxdt
\leq C_2\,.\eqno(3.5)$$
Call $t\mapsto x(t)\doteq\Phi^\nu_t(y)$ the solution of
$$\dot x=f_\nu(t,x),\qquad\qquad x(0)=y\,.$$
Moreover, assume that the fluxes $\Phi_t^\nu$ are all
nearly incompressible, so that, for every bounded set $A\subset\R^m$,
$${1\over C_3}\,\hbox{meas}(A)~\leq~\hbox{meas}\big(\Phi_t^\nu(A)\big)~\leq~
C_3\,\hbox{meas}(A)\,,\eqno(3.6)$$
for some constant $C_3$ and all $t\in [0,T]$, $\nu\geq 1$. 
\v
\n{\bf Conjecture:}~~By possibly exstracting a subsequence,
one has the convergence
$$\Phi^\nu\to\Phi\qquad\hbox{in}~~~\L^1_{loc}$$
for some measurable flux $\Phi$, also satisfying (3.6).

Notice that, by a well known compactness theorem, for some $f\in\L^\infty$
we certainly have
$$f_\nu\to f\qquad\hbox{in}~~~\L^1_{loc}\,.$$
The existence of a limit flux $\Phi$ would provide a new existence
result concerning
the Cauchy problem for the discontinuous O.D.E.
$$\dot x=f(t,x)\,,\qquad\qquad x(0)=y\eqno(3.7)$$
valid for a.e.~initial data $y\in\R^m$, somewhat extending
the result in [3]. 
\vsk
\c{\medbf References}
\v
\i{[1]} A.Bressan, {\it Hyperbolic Systems of Conservation Laws.
The One Dimensional Cauchy Problem}, Oxford University Press, 2000.
\v
\i{[2]} C.~Dafermos, {\it Hyperbolic Conservation Laws in Continuum
Physics}, Springer-Verlag, Berlin 1999.
\v
\i{[3]} R.~DiPerna and P.~L.~Lions, Ordinary differential equations,
transport theory and Sobolev spaces, {\it Invent. Math.} {\bf 98} 
(1989), 511-517.
\v
\i{[4]} S.~Kruzhkov,
First-order quasilinear equations with several space
variables, {\it Math. USSR Sbornik} {\bf 10} (1970), 217--273.
\v
\i{[5]} D.~Serre, {\it Systems of Conservation Laws I, II},
Cambridge University Press, 2000.

\bye